\documentclass[11pt]{article}
\usepackage{amsmath}
\usepackage{amsfonts}
\usepackage{amssymb}
\usepackage{amsthm}
\usepackage{graphicx}

\providecommand{\U}[1]{\protect\rule{.1in}{.1in}}
\newtheorem{theorem}{Theorem}[section]

\newtheorem{proposition}[theorem]{Proposition}
\theoremstyle{definition}
\newtheorem{definition}[theorem]{Definition}
\theoremstyle{remark}
\newtheorem{example}[theorem]{Example}
\newtheorem{remark}[theorem]{Remark}
\numberwithin{equation}{section}
\newenvironment{Proof}[1][Proof]{\noindent\textbf{#1.} }{\ \rule{0.5em}{0.5em}}
\begin{document}

\title{Analytic approximation of solutions of parabolic partial differential
equations with variable coefficients}
\author{Vladislav V. Kravchenko, Josafath A. Otero, Sergii M. Torba\\
{\small Departamento de Matem\'{a}ticas, CINVESTAV del IPN, Unidad Quer\'{e}%
taro, }\\
{\small Libramiento Norponiente No. 2000, Fracc. Real de Juriquilla, Quer%
\'{e}taro, Qro. C.P. 76230 MEXICO}\\
{\small e-mail: vkravchenko@math.cinvestav.edu.mx,
josafath@math.cinvestav.edu.mx, storba@math.cinvestav.edu.mx}}
\maketitle

\begin{abstract}
A complete family of solutions for the one-dimensional reaction-diffusion
equation
\begin{equation}
u_{xx}(x,t)-q(x)u(x,t) = u_t(x,t) \label{abs1}
\end{equation}%
with a coefficient $q$ depending on $x$ is constructed. The solutions
represent the images of the heat polynomials under the action of a
transmutation operator. Their use allows one to obtain an
explicit solution of the noncharacteristic Cauchy problem for
equation \eqref{abs1} with sufficiently regular Cauchy data as well as to solve
numerically initial boundary value problems. In the paper the Dirichlet boundary conditions are considered however the proposed method can be easily extended onto other standard boundary conditions. 
The proposed numerical method
is shown to reveal good accuracy.
\end{abstract}

\section{Introduction}

In the present work a complete system of solutions of a one-dimensional
reaction-diffusion equation with a variable coefficient%
\begin{equation}
u_{xx}(x,t)-q(x)u(x,t) = u_t(x,t) \label{eq rd}
\end{equation}%
considered on $\Omega :=(-b,b)\times (0,\tau )$ is obtained. We assume that the potential $q\in C[-b,b]$ may be complex valued. The completeness of the system is with respect to the uniform norm in the closed rectangle $\bar \Omega$. The system of solutions is shown to be useful for
uniform approximation of solutions of initial boundary value problems for (\ref{eq rd})\footnote{In the paper the Dirichlet boundary conditions are considered however the proposed method can be easily extended onto other standard boundary conditions.} as well as for explicit solution of the noncharacteristic
Cauchy problem (see \cite{cannon}) for (\ref{eq rd}) in terms of the
formal powers arising in the spectral parameter power series
(SPPS) method (see \cite{Kr2008}, \cite{KrP2010SPPS}, \cite{KrTr2014SPPSM}).

The complete system of solutions is constructed with the aid of the
transmutation operators relating (\ref{eq rd}) with the
heat equation (see e.g., \cite{colton}, \cite{marchenko}, \cite{KKrTrT2013}%
). The possibility to construct complete systems of solutions by means of
transmutation operators was proposed and explored in \cite{colton}, and the
approach developed in \cite{colton} requires the knowledge of the
transmutation operators. In the present work using a mapping property of the
transmutation operators discovered in \cite{CKrTr2012} we show that the
construction of the complete systems of solutions for equations of the form (%
\ref{eq rd}), representing transmuted heat polynomials, can be realized with
no previous construction of the transmutation operator. Moreover, the use of
the mapping property leads to an explicit solution of the noncharacteristic
Cauchy problem for (\ref{eq rd}) with Cauchy data belonging to a Holmgren
class \cite{cannon}.

We illustrate the implementation of the complete system of the transmuted
heat polynomials by a numerical solution of an initial boundary value
problem for (\ref{eq rd}). The approximate solution is sought in the form of
a linear combination of the transmuted heat polynomials and the initial and
boundary conditions are satisfied by a collocation method. A remarkable
accuracy is achieved in few seconds using Matlab 2012 on a usual PC.

Besides this Introduction the paper contains five sections. In Section \ref{Sect2} we
recall the transmutation operators and some of their properties. In Section
\ref{Sect3} an explicit solution of the noncharacteristic Cauchy problem for (\ref{eq
rd}) is obtained. In Section \ref{Sect4} we prove the completeness of the transmuted
heat polynomials. In Section \ref{Sect5} the numerical method for solving initial
boundary value problems for (\ref{eq rd}) implementing the transmuted heat
polynomials is discussed. Section \ref{Sect6} presents a numerical illustration.

\section{Transmutation operators and formal powers}\label{Sect2}

\subsection{System of recurrent integrals}

Let $q$ be a continuous complex valued function defined on the segment $%
[-b,b]$. Throughout the paper we suppose that $f$ is a nonvanishing solution
of the equation
\begin{equation}
f^{\prime \prime }-q(x)f=0  \label{particular_solution}
\end{equation}%
on $(-b,b)$ such that $f(0)=1$ and $f^{\prime }(0)=\alpha $ where $\alpha $
is a complex number. In \cite{KrP2010SPPS} the existence of such
solution was proved.

Consider two sequences of recurrent integrals (see \cite{Kr2008}, \cite%
{KrP2010SPPS}, \cite{KrMT2012})
\begin{align*}
X^{(0)} \equiv 1,\quad
X^{(n)}(x)&=n\int_{0}^{x}X^{(n-1)}(s)(f^{2}(s))^{(-1)^{n}}ds,\qquad x\in[-b,b],\ n\in
\mathbb{N},\\  
\tilde{X}^{(0)}\equiv 1,\quad \tilde{X}^{(n)}(x)&=n\int_{0}^{x}\tilde{X}%
^{(n-1)}(s)(f^{2}(s))^{(-1)^{n-1}}ds,\qquad x\in[-b,b],\ n\in
\mathbb{N}.  
\end{align*}

\begin{definition}
The family of functions $\{\varphi _{k}\}_{k=0}^{\infty }$ constructed
according to the rule
\begin{equation}
\varphi _{k}(x)=
\begin{cases}
f(x)X^{(k)}(x), & k\text{ odd} \\
f(x)\tilde{X}^{(k)}(x), & k\text{ even}%
\end{cases}\label{formal powers}
\end{equation}%
is called the system of formal powers associated with $f$.
\end{definition}

The formal powers arise in the spectral parameter power series (SPPS)
representation for solutions of the one-dimensional Schr\"odinger equation
(see \cite{Kr2008}, \cite{KrP2010SPPS}).

\subsection{The transmutation operator}

For any $q\in C[-b,b]$ it is a well known result \cite[Chapter 1]{marchenko}
that there exists a function $(x,s) \mapsto K(x,s)$ defined on the domain $0\leq
\left\vert s\right\vert \leq |x|\leq b$, continuously differentiable, such that the equality
\begin{equation*}
ATv=TBv  
\end{equation*}
is valid for all $v\in C^2[-b,b]$, where $A:=\frac{\partial ^{2}}{\partial
x^{2}}-q$, $B:=\frac{\partial ^{2}}{\partial x^{2}}$ and $T$ has the form of
a second kind Volterra integral operator 
\begin{equation*}
Tv(x):=v(x)+\int_{-x}^{x}K(x,s)v(s)ds.
\end{equation*}%
The operator $T$ is called \textbf{transmutation operator}.
Moreover, the function $K$ is not unique and can be chosen so that $T[1]=f$
(see, e.g., \cite{KT2016BSMM}). When $q\in C^{1}[-b,b]$ such function $K$ is
the unique solution of the Goursat problem
\begin{gather*}
K_{xx}(x,s)-q(x)K(x,s) =K_{ss}(x,s), \\
K(x,x)=\frac{\alpha }{2}+\frac{1}{2}\int_{0}^{x}q(y)dy, \qquad x\in \lbrack -b,b],
\\
K(x,-x)=\frac{\alpha }{2}, \qquad x\in \lbrack -b,b].
\end{gather*}%
For any $q\in C[-b,b]$ the kernel $K$ can be defined as $K(x,s)=H\big(\frac{%
x+s}{2},\frac{x-s}{2}\big)$, $|s|\leq |x|\leq b$, with $H$ being the unique
solution of the Goursat problem
\begin{gather*}
H_{uv}(u,v)=q(u+v)H(u,v),\\
H(u,0)=\frac{\alpha }{2}+\frac{1}{2}\int_{0}^{u}q(s)\,ds,\qquad H(0,v)=\frac{%
\alpha }{2}.  
\end{gather*}%
If the potential $q$ is $n$ times continuously differentiable on $(-b,b)$, the kernel $%
K(x,s)$ is $n+1$ times continuously differentiable with respect to both
independent variables.

The following mapping property of the operator $T$ is used throughout the
paper.

\begin{proposition}[\cite{CKrTr2012}]
\begin{equation}
T[x^{k}]=\varphi _{k}(x),\quad \forall \;k\in \mathbb{N}_{0}.
\label{mapping propertie}
\end{equation}
\end{proposition}

The inverse operator $T^{-1}$ also has the form of a second kind Volterra integral operator and satisfy the following correspondence of the initial values, see \cite{KrTr2012D}
\begin{equation}
v(0) = u(0),\qquad v'(0)=u'(0)-\alpha u(0), \label{initial values}
\end{equation}
where $v:=T^{-1}u$.

\section{The noncharacteristic Cauchy problem for (\protect\ref{eq rd})}\label{Sect3}

In this section an explicit solution of the noncharacteristic Cauchy problem
for (\ref{eq rd}) in terms of the formal powers $\varphi _{k}$ is obtained.
This fact is a direct consequence of the mapping property (\ref{mapping
propertie}).

\begin{definition}[\protect\cite{cannon}]
For the positive constants $\gamma _{1}$, $\gamma _{2}$ and $C_{1}$, the
Holmgren class $H(\gamma _{1},\gamma _{2},C_{1},t_{0})$ is the set of
infinitely differentiable functions $v$ defined on $\left\vert
t-t_{0}\right\vert <\gamma _{2}$ that satisfy
\begin{equation*}
\left\vert v^{(j)}(t)\right\vert \leq C_{1}\gamma _{1}^{-2j}(2j)!,\quad
j=0,1,\ldots
\end{equation*}%
for all $t\in \left\vert t-t_{0}\right\vert <\gamma _{2}$.
\end{definition}

\begin{proposition}
Let $q\in C[-b,b]$ and $u(x,t)$ be a solution of the noncharacteristic
Cauchy problem
\begin{align}
u_{xx}(x,t)-q(x)u(x,t) &=u_{t}(x,t),\qquad -b<x<b,\quad \left\vert t\right\vert <\tau
\label{eq cauchy} \\
u(0,t) &=F(t),\qquad \left\vert t\right\vert <\tau \\
u_{x}(0,t) &=G(t),\qquad \left\vert t\right\vert <\tau
\end{align}%
where $F,G\in H(b,\tau ,C,0)$, $C>0$. Then the series
\begin{equation*}
\sum_{j=0}^{\infty }\left[\frac{F^{(j)}(t)}{(2j)!}\left(\varphi _{2j}(x)-\frac{\alpha}{2j+1}\varphi_{2j+1}(x)\right)+\frac{G^{(j)}(t)%
}{(2j+1)!}\varphi _{2j+1}(x)\right]
\end{equation*}%
converges uniformly and absolutely for $\left\vert x\right\vert \leq r<b$ to
the solution $u(x,t)$ where $\varphi _{k}$ are the formal powers
\eqref{formal powers}.
\end{proposition}

\begin{remark}
Note that the functions $\varphi _{2j}(x)-\frac{\alpha}{2j+1}\varphi_{2j+1}(x)$ coincide with the formal powers \eqref{formal powers} constructed starting with the particular solution $g$ of \eqref{particular_solution} satisfying the initial conditions $g(0)=1$ and $g'(0)=0$, see \cite[Proposition 4.7]{KrTr2013}.
\end{remark}

\begin{Proof}
Let $u(x,t)$ be a solution of ($\ref{eq cauchy})$. Consider the function
$h:=T^{-1} u$, where the operator $T^{-1}$ is applied with respect to the variable $x$. The function $h$ is a solution of the heat equation, c.f., \cite[Theorem 2.1.2]{colton} and due to \eqref{initial values} satisfies the following noncharacteristic Cauchy problem
\begin{align*}
h_{xx} &=h_{t},\qquad -b<x<b,\quad \left\vert t\right\vert <\tau \\
h(0,t) &=F(t),\qquad \left\vert t\right\vert <\tau \\
h_{x}(0,t) &=G(t)-\alpha F(t),\qquad \left\vert t\right\vert <\tau.
\end{align*}
Since $G-\alpha F\in H(b,\tau ,(1+\alpha)C,0)$, the solution of this problem is given by the absolutely and uniformly convergent series for $\left\vert
x\right\vert \leq r<b$ (see, e.g., \cite{cannon})
\begin{equation*}
h(x,t)=\sum_{k=0}^{\infty }\left(\frac{F^{(k)}(t)}{(2k)!}x^{2k}+\frac{G^{(k)}(t)-\alpha F^{(k)}(t)}{%
(2k+1)!}x^{2k+1}\right).
\end{equation*}%
Due to (\ref{mapping propertie}) we obtain
\begin{equation*}
u(x,t)=Th(x,t)=\sum_{k=0}^{\infty }\left[\frac{F^{(k)}(t)}{(2k)!}\left(\varphi _{2k}(x)-\frac{\alpha}{2k+1}\varphi_{2k+1}(x)\right)+%
\frac{G^{(k)}(t)}{(2k+1)!}\varphi _{2k+1}(x)\right].
\end{equation*}%
This series converges uniformly and absolutely for $\left\vert x\right\vert
\leq r<b$ due to the uniform boundedness of $T$ and of its inverse.
\end{Proof}

\section{Transmuted heat polynomials}\label{Sect4}

In this section a complete system of solutions of (\ref{eq rd}) is presented.

Consider the heat polynomials (see, e.g., \cite{WidBl1959}) defined by
\begin{equation}
h_{n}(x,t)=n!\sum_{k=0}^{\left[ \frac{n}{2}\right] }\frac{t^{k}x^{n-2k}}{%
k!(n-2k)!},\quad n\in \mathbb{N}_{0}.  \label{heat polynomials}
\end{equation}%
Due to (\ref{mapping propertie}) we obtain that the functions
\begin{equation}
u_{n}(x,t)=n!\sum_{k=0}^{\left[ \frac{n}{2}\right] }\dfrac{t^{k}\varphi
_{n-2k}(x)}{k!(n-2k)!},\quad n\in \mathbb{N}_{0}
\label{transformed heat polynomials}
\end{equation}%
are solutions of (\ref{eq rd}) for all $-b<x<b$ and $t>0$. Indeed, we have
that $u_{n}=Th_{n},$ $n\in \mathbb{N}_{0}$ and
\begin{equation*}
\frac{\partial ^{2}u_{n}}{\partial x^{2}}-q(x)u_{n}(x,t)-\frac{\partial u_{n}%
}{\partial t}(x,t)=T\left( \frac{\partial ^{2}h_{n}}{\partial x^{2}}-\frac{%
\partial h_{n}}{\partial t}(x,t)\right) =0.
\end{equation*}

The completeness of the system of the heat polynomials with respect to the maximum norm proved in \cite{colton} and the
uniform boundedness of $T$ and $T^{-1}$ imply the completeness of (\ref%
{transformed heat polynomials}) in the space of classical solutions of (\ref%
{eq rd}). Thus, the following statement is true.

\begin{theorem}
\label{thm thp} Let $u(x,t)$ be continuous in $\bar{\Omega}$ and satisfy 
\eqref{eq rd} in $\Omega $. Then given $\varepsilon >0$ there exists $N\in
\mathbb{N}$ and constants $a_{0},a_{1},\ldots ,a_{N}$ such that
\begin{equation*}
\max_{\bar{\Omega}}\left\vert
u(x,t)-\sum_{n=0}^{N}a_{n}u_{n}(x,t)\right\vert <\varepsilon .
\end{equation*}
\end{theorem}

\begin{Proof}
Choose $\varepsilon >0$. Consider $h(x,t)=T^{-1}u(x,t)$. Due to the
completeness of (\ref{heat polynomials}), for any $\varepsilon _{1}>0$ there
exists $N\in \mathbb{N}$ and constants $a_{0},a_{1},\ldots ,a_{N}$ such that
\begin{equation*}
\max_{\bar{\Omega}}\left\vert
h(x,t)-\sum_{n=0}^{N}a_{n}h_{n}(x,t)\right\vert <\varepsilon _{1}.
\end{equation*}%
Then
\begin{equation*}
\max_{\bar{\Omega}}\left\vert
u(x,t)-\sum_{n=0}^{N}a_{n}u_{n}(x,t)\right\vert =\max_{\bar{\Omega}%
}\left\vert Th(x,t)-\sum_{n=0}^{N}a_{n}Th_{n}(x,t)\right\vert \leq
C\varepsilon _{1}
\end{equation*}%
where the constant $C$ is the uniform norm of $T$. The choice of $%
\varepsilon _{1}=\varepsilon /C$ finishes the proof.
\end{Proof}

\section{Solution of initial boundary value problems for (\protect\ref{eq rd}%
)}\label{Sect5}

Consider the problem to find the solution of the equation
\begin{equation}
u_{xx}(x,t)-q(x)u(x,t)=u_t(x,t),\quad (x,t)\in \Omega  \label{eq Initial-Dirichlet}
\end{equation}%
subject to the Dirichlet boundary conditions
\begin{equation}
u(-b,t)=\psi _{1}(t),\quad u(b,t)=\psi _{2}(t),\quad t\in \lbrack 0,\tau ]
\label{dirichlet condition}
\end{equation}%
and the initial condition
\begin{equation}
u(x,0)=\varphi (x),\quad x\in \lbrack -b,b]  \label{initial condition}
\end{equation}%
where $\psi _{1}$, $\psi _{2}$ and $\varphi $ are continuously
differentiable functions satisfying the compatibility conditions
\begin{equation*}
\psi _{1}(0)=\varphi (-b),\quad \psi _{2}(0)=\varphi (b).
\end{equation*}%
The problem \eqref{eq Initial-Dirichlet}--\eqref{initial condition} possesses a
unique solution and depends continuously on the data (see, e.g., \cite%
{vladimirov}).

The result of Theorem \ref{thm thp} suggests the following simple method to
approximate the solution of problem (\ref{eq Initial-Dirichlet})--(\ref%
{initial condition}). The approximate solution $\tilde{u}$ is sought in the
form
\begin{equation}
\tilde{u}(x,t)=\sum_{n=0}^{N}a_{n}u_{n}(x,t).  \label{approximated solution}
\end{equation}%
Since every $u_{n}$ is a solution of (\ref{eq Initial-Dirichlet}), their
linear combination satisfies (\ref{eq Initial-Dirichlet}) as well. The
coefficients $\left\{ a_{n}\right\} _{n=0}^{N}$ are sought in such way that $%
\tilde{u}$ satisfy the initial and the boundary conditions approximately.
For this we used the collocation method. $M$ points $\{(x_{i},t_{i})%
\}_{i=1}^{M}$ are chosen on the parabolic boundary $(\{-b\}\times \lbrack
0,\tau ])\cup ([-b,b]\times \{0\})\cup (\{b\}\times \lbrack 0,\tau ])$ in
order to construct a linear system of equations for the coefficients $%
\left\{ a_{n}\right\} _{n=0}^{N}$
\begin{equation}
\sum_{n=0}^{N}a_{n}u_{n}(x_{i},t_{i})=
\begin{cases}
\psi _{1}(t_{i}), & x_{i}=-b \\
\varphi (x_{i}), & t_{i}=0 \\
\psi _{2}(t_{i}), & x_{i}=b
\end{cases}
,\qquad i=1,\ldots ,M.  \label{coefficients}
\end{equation}%
The system (\ref{coefficients}) is the result of imposing the conditions (%
\ref{dirichlet condition}), (\ref{initial condition}) onto the approximate
solution (\ref{approximated solution}). Using the pseudoinverse matrix the
system (\ref{coefficients}) is solved, and the approximate solution (\ref%
{approximated solution}) is computed on $\Omega $ using the obtained
coefficients $\left\{ a_{n}\right\} _{n=0}^{N}$ and the definition of the
transmuted heat polynomials (\ref{transformed heat polynomials}).

Needless to add that the same approach is applicable to other kinds of
boundary conditions.

\section{Numerical illustration}\label{Sect6}

We present a numerical example of the application of the method described in
the previous section. It reveals a remarkable accuracy with very little
computational efforts. The implementation was realized in Matlab 2012.

On the first step a nonvanishing solution $f$ of \eqref{particular_solution}
was computed using the SPPS method (see \cite{KrP2010SPPS}, \cite%
{KrTr2014SPPSM}). The formal powers $\varphi _{k}$ were constructed like in
\cite{KKrTrT2013} using the \texttt{spapi} and \texttt{fnint} Matlab
routines from the spline Toolbox. Then the transmuted heat polynomials \eqref{transformed heat polynomials} were calculated. In order to obtain a unique
solution of the system \eqref{coefficients} $N+1$ equally spaced points on the
parabolic boundary were chosen. Finally, the approximate solution (\ref%
{approximated solution}) was computed on a mesh of $200\times 100$ points in the
interior of the rectangle and compared with the corresponding exact solution.

\begin{example}\label{Ex1}
Consider the initial Dirichlet problem
\begin{align}
u_{xx}(x,t)-x^{2}u(x,t) &=u_{t}(x,t),\qquad (x,t)\in (-1,1)\times (0,1),  \label{Example_eq}
\\
u(x,0) &=e^{-0.5x^{2}},\qquad x\in \lbrack -1,1],  \label{Example_initial} \\
u(-1,t) &=u(1,t)=e^{-0.5-t},\qquad t\in \lbrack 0,1].
\label{Example_frontera}
\end{align}

\begin{figure}[tbh]
\centering
\includegraphics[height=3.18in, width=4.76in]{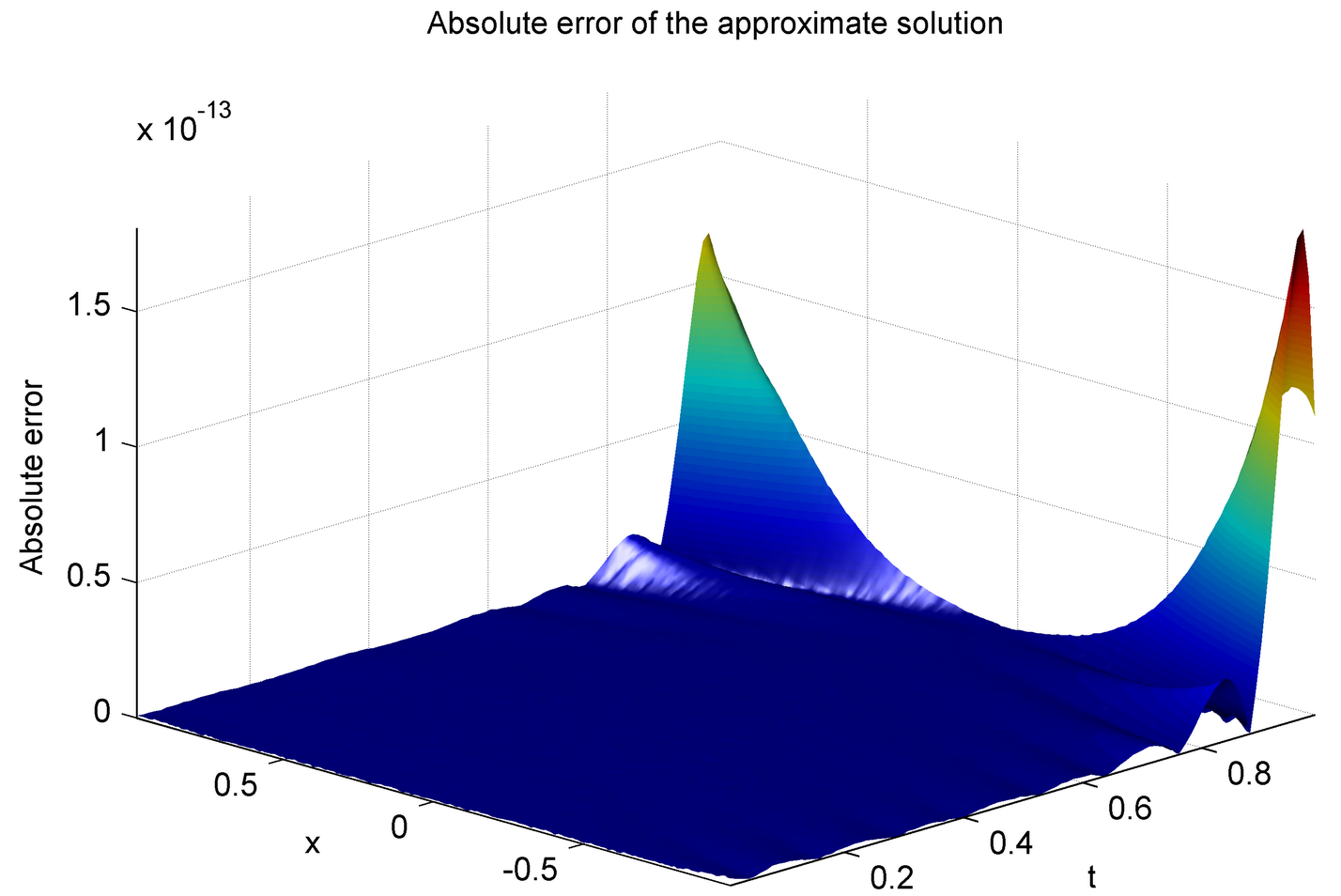}
\caption{The absolute value of the difference $| u(x,t)-u_{26}(x,t)|$  between the exact and the approximate solutions for the problem \eqref{Example_eq}--\eqref{Example_frontera}.}
\label{Fig abserror}
\end{figure}

\begin{table}[tbh]
\centering
\begin{tabular}{|c|c|c|c|}
\hline
$N$ & Max.\ absolute error & Max.\ relative error & Cond.\ number\\ \hline
5 & $2.3\cdot 10^{-2}$ & $6.2\cdot 10^{-2}$ & $55.7$\\ \hline
10 & $1.4\cdot 10^{-4}$ & $5.5\cdot 10^{-4}$ & $1.82\cdot 10^{5}$\\ \hline
15 & $9.6\cdot 10^{-8}$ & $4.0\cdot 10^{-7}$ & $3.65\cdot 10^{9}$\\ \hline
20 & $2.0\cdot 10^{-10}$ & $8.6\cdot 10^{-10}$ & $1.76\cdot 10^{14}$\\ \hline
23 & $7.6\cdot 10^{-13}$ & $3.2\cdot 10^{-12}$ & $1.67\cdot 10^{17}$\\ \hline
26 & $1.8\cdot 10^{-13}$ & $7.8\cdot 10^{-13}$ & $2.59\cdot 10^{23}$\\ \hline
29 & $2.5\cdot 10^{-12}$ & $1.1\cdot 10^{-11}$ & $1.23\cdot 10^{23}$\\ \hline
34 & $1.7\cdot 10^{-10}$ & $7.3\cdot 10^{-10}$ & $1.37\cdot 10^{25}$\\ \hline
39 & $2.3\cdot 10^{-9}$ & $9.8\cdot 10^{-9}$ & $2.86\cdot 10^{29}$\\ \hline
50 & $4.7\cdot 10^{-10}$ & $2.1\cdot 10^{-9}$ & $1.25\cdot 10^{41}$\\ \hline
75 & $6.1\cdot 10^{-11}$ & $2.7\cdot 10^{-10}$ & $2.89\cdot 10^{73}$\\ \hline
100 & $2.8\cdot 10^{-10}$ & $1.2\cdot 10^{-9}$ & $5.63\cdot 10^{105}$\\ \hline
\end{tabular}
\caption{Maximal absolute and relative errors of the approximate solution and condition number of the matrix in \eqref{coefficients} for the problem \eqref{Example_eq}--\eqref{Example_frontera} obtained for different values of $N$ in \eqref{approximated solution}.}
\label{Tab Errors}
\end{table}

The exact solution of this problem has the form
\begin{equation*}
u(x,t)=\exp \left( -\frac{1}{2}x^{2}-t\right) .
\end{equation*}

The distribution of the absolute error of the approximate solution for $N=26$ is presented on Figure \ref{Fig abserror}. The maximum absolute error of
the approximate solution is of order $10^{-13}$.

It is often stated that boundary collocation methods (in particular, the heat polynomials method) lead to ill-conditioned systems of linear equations, see \cite{ChH1994}, \cite{KZ2009}, \cite{LLHCh2008}. It is also the case for the proposed method. As is illustrated in Table \ref{Tab Errors}, the condition number of the matrix in \eqref{coefficients} grows rather fast. Nevertheless, the straightforward implementation of the proposed method presented no numerical difficulties.
The convergence and the robustness of the method are illustrated in Table \ref{Tab Errors} where the maximum absolute and the maximum relative error of the approximate solution for different values of $N$ used for approximation \eqref{approximated solution} are presented. As one can appreciate, the convergence rate for small values of $N$ is exponential. And 
even taking values of $N$ much larger than the optimum one do not lead to any problem for collocation method nor to significant precision lost. Moreover, a simple test based on the accuracy of fulfilment of the initial and boundary conditions \eqref{Example_initial}--\eqref{Example_frontera} can be utilized to estimate both the optimal $N$ and the accuracy of the obtained approximate solution.
\end{example}

\section*{Conclusions}
A complete system of solutions of equation (\ref{eq rd}) is obtained. The
solutions represent the images of the heat polynomials under the action of the
transmutation operator. They are shown to be convenient for uniform
approximation of solutions of initial boundary value problems for
(\ref{eq rd}) as well as for explicit solution of the noncharacteristic Cauchy
problem. Besides the Dirichlet boundary conditions considered in this paper
the method is applicable to other standard boundary conditions. The complete
system of solutions obtained can be used for solving moving and free boundary
problems \cite{KKT2017}.

\section*{Conflicts of Interest}
The authors declare that there is no conflict of interest regarding the publication of this paper.

\section*{Acknowledgments}
Research was supported by CONACYT, Mexico via the project 222478.

\end{document}